# Ramsey Theory and Geometry of Closed Loops

Nir Shvalb[1], Mark Frenkel[2], Shraga Shoval[3], Edward Bormashenko[2*]

[1]Department of Mechanical Engineering & Mechatronics, Faculty of Engineering, Ariel University, P.O.B.3, Ariel, 407000, Israel

[2]Chemical Engineering Department, Engineering Faculty, Ariel University, P.O.B. 3, 407000 Ariel, Israel

[3]Department of Industrial Engineering and Management, Faculty of Engineering, Ariel University, P.O.B. 3, Ariel 407000, Israel

*Correspondence: edward@ariel.ac.il

**Abstract**

We apply the Ramsey theory to the analysis of geometrical properties of closed contours. Consider a set of six points placed on a closed contour. The straight lines connecting these points are $y_{ik}(x) = \alpha_{ik}x + \beta_{ik}$ $(i, k = 1\dots 6), \alpha_{ik} \neq 0$. We color the edges connecting the points for which $\alpha_{ik} > 0$ holds with red, and the edges for which $\alpha_{ik} < 0$ with green with red. At least one monochromic triangle should necessarily appear within the curve (according to the Ramsey number $R(3,3) = 6$). This result is immediately applicable for the analysis of dynamical billiards. The second theorem emerges from the combination of the Jordan curve and Ramsey theorem. The closed curve is considered. We connect the points located within the same region with green links and the points placed within the different regions with red links. In this case, transitivity/intransitivity of the relations between the points should be considered. Ramsey constructions arising from the differential geometry of closed contours are discussed.

1. Introduction

Ramsey theory is a branch of graph theory that focuses on the appearance of interconnected substructures within a structure/graph of a known size [1-9]. Ramsey theory states that any structure will necessarily contain an interconnected substructure [1-6]. Ramsey's theorem, in one of its graph-theoretic forms, states that one will find monochromatic cliques in any edge labelling (e.g. with colors) of a sufficiently large complete graph [6]. One more example of the Ramsey-like thinking is delivered by the van der Waerden's theorem: colorings of the integers by finitely many colors must have long monochromatic arithmetic progressions [6]. An accessible introduction to the

Ramsey theory is found in refs. 1-8. More rigorous approach is laid out in refs. 7-8. Problems in Ramsey theory typically ask a question of the form: "how big must some structure be to guarantee that a particular property holds?" We apply the Ramsey theory to analysis of the geometrical properties of closed contours; thus, bridging between the discrete and continuous mathematic approaches.

## 2. Results

### 2.1. Ramsey theory and geometry of closed curves

Consider the closed curve depicted in **Figure 1** with a solid black line. Let us choose six points placed on the curve and number them $i = 1 \ldots 6$ as shown in **Figure 1**.

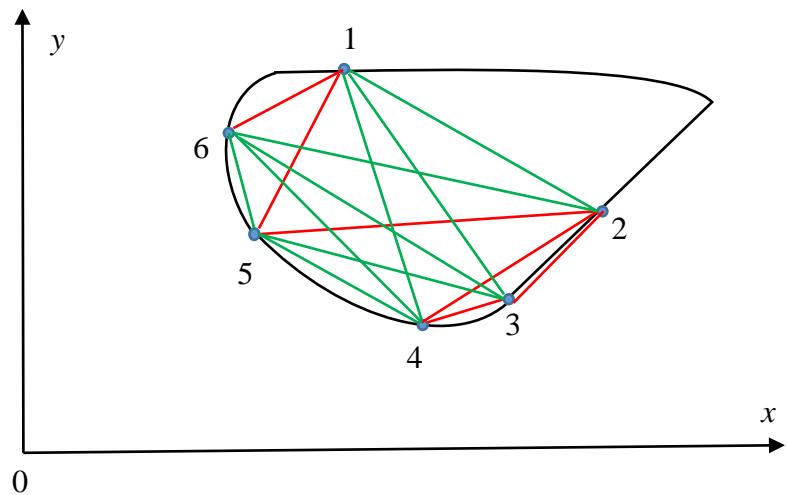

**Figure 1**. Ramsey construction for six points placed on the closed curve.

We connect each pair of points in pairs with the straight colored segments, as shown in **Figure 1**. The equations of these straight lines are given by Eq. 1:

$$y_{ik}(x) = \alpha_{ik}x + \beta_{ik} \ (i, k = 1\ldots 6), \tag{1}$$

where $y_{ik}(x)$ is the equation of the straight line connecting points numbered $i$ and $k$ correspondingly; $\alpha_{ik}$ is the slope of the straight line and $\beta_{ik}$ is its $y$-intercept. It is always possible to choose six points laying in a closed 2D curve in such a way that $\alpha_{ik} \neq 0$ holds.

**Proof:** Jordan's theorem holds for any closed simple curves in the plane including exotic ones like the Warsaw circle. Nevertheless, for simplicity let restrict ourselves to well behaved piecewise-polynomial curves F over an algebraically closed field. Since the curve {y=0} and F do not share a common factor F ∩ {y=0} is a finite set of points (a generalized case is provided in ref. 13). Thus, for any piecewise-polynomial curve excluding a finite number of points will do for the above statement to hold.

Thus, two kinds of slope, namely $\alpha_{ik} > 0$ and $\alpha_{ik} < 0$ are possible ($\alpha_{ik} \neq 0$ is assumed). Lets color the edges connecting the points, numbered $i$ an k for which $\alpha_{ik} > 0$ holds with red, and the edges for which $\alpha_{ik} < 0$ with green color, as shown in **Figure 1**. Thus, the complete bi-color graph is created, and according to the Ramsey theorem, at least one (red or green) monochromic triangle should necessarily appear within the graph; due to the fact that the Ramsey number is $R(3,3) = 6$. Indeed, the triangle "456" appearing in **Figure 1** is built of green edges only. And this will be true for any closed contour.

### 2.2. Ramsey theory and dynamical billiards

The aforementioned theorem enables re-shaping of the problem of dynamical billiards in a spirit of the Ramsey theory. Consider the dynamical billiard depicted in **Figure 2**.

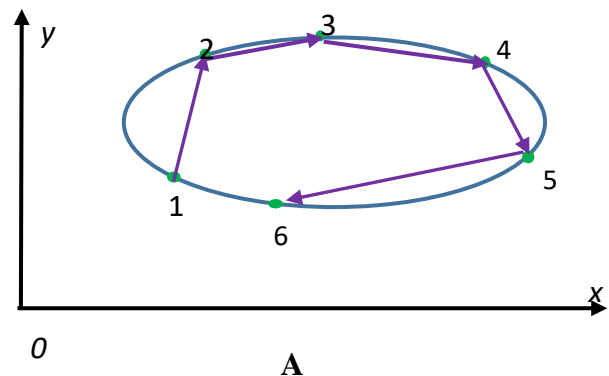

A

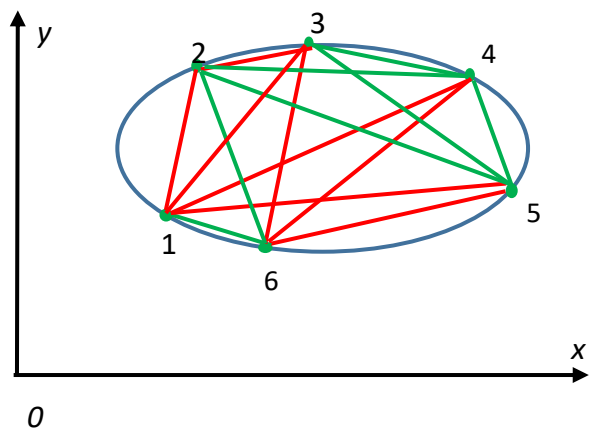

**B**

**Figure 2**. Dynamic billiard as the Ramsey system. **A**. Trajectory of the point within the pool (billiard) is shown with violet arrows. **B**. Complete graph arising from six points, depicting reflections from the boundary.

In the dynamical billiard, a particle moves along a straight line and is reflected from the boundaries. Billiards are Hamiltonian idealizations of known billiard game, in which the boundaries have general geometric shape (rather than the rectangular shape). In **Figure 2A** a particle alternates between free motion (presupposed to be a straight line, depicted with violet arrows) and specular reflections from a boundary. This class of systems is called a dynamical billiard [10]. Dynamical billiards are described by mathematical models which appear in a diversity of physical phenomena. The dynamical properties of such models are determined by the shape of the walls of the container, and they may vary from completely regular (integrable) to fully chaotic [10-12]. Thus, consider the simplest system in which the particles moves in a 2D container and collide with its walls/boundary, as shown in **Figure 2A**. The reflection points are shown with green circles. The boundary is supposed to be a closed curve, shown with a blue sold line. Consider the set of first six reflections from the boundary shown with green circles in **Figure 2A** and numbered $i = 1 \ldots 6$. We connect the points in pairs with straight line segments, as shown in **Figure 2B**. The equations of these straight lines again are provided by Eq. 1. Following the Ramsey procedure and approach discussed in Section 2.1 we color the edges connecting the reflection points, numbered $i$ an $k$ for which $\alpha_{ik} > 0$ holds with red, and the links for which $\alpha_{ik} < 0$ with green color, as shown in **Figure 2B**. Again, $\alpha_{ik} \neq 0$ is presumed. Thus, the complete bi-color graph emerges, and according to the Ramsey theorem, at least one (red or green) monochromic triangle should necessarily appear within the graph (the Ramsey number $R(3,3) = 6$). Indeed, triangle "123" appearing in **Figure 2B** is built of red edges only, and triangles "245" and "345" are built from the green edges only, i.e. two monochromatic triangles are recognized in **Figure 2B**. This result will be true for any closed boundary, for any starting point of the body moving within the billiard. Moreover, the reflections may be not exactly specular. We also do not specify Hamiltonian, representing the energy of the particle. Any set of six reflection points located on the closed boundary will generate the complete graph, which fulfils the

Ramsey theorem, and at least one monochromatic triangle will necessary appear in the graph.

Dynamic billiards are of much interest in a view of deterministic chaos, revealed for these mechanical systems [10-12]. Ramsey theorem states that this chaos will never be complete, substructures built of the segments of the trajectory of the particle will necessary appear, as illustrated in **Figure 2B**.

*2.3. Transitive Ramsey numbers and differential geometry of closed curves*

Consider now the closed curve, shown in **Figure 3**. The curvature of this curve is positive at the black segment of the curve and negative at the violet segment of the curve. Recall that the signed curvature $\kappa$ is given by Eq. 2:

$$\kappa = \frac{y''}{\left(1+y'^2\right)^{3/2}} \qquad (2)$$

Consider a set of *n* points placed on the curve. Some of the points are placed on the segments with the positive curvature (they are denoted with symbol "+" in **Figure 3**), and some are placed on the segments with negative curvature (they are denoted with symbol "-" in **Figure 3**). According to the Ramsey approach, we connect with the green edge the points, with the same sign; and differently signed points we connect with the red edge, as shown in **Figure 3**. This procedure gives rise to the complete graph. Let us answer the following question: what it the minimal number of points which necessarily provides formation of the monochromatic triangle in the aforementioned complete graph? It also seems from first glance, that the answer is supplied by the Ramsey theorem and $R(3,3) = 6$ is kept. However, this answer is wrong, due to the fact, that the presented relations between the vertices of the graph are transitive/intransitive in their nature.

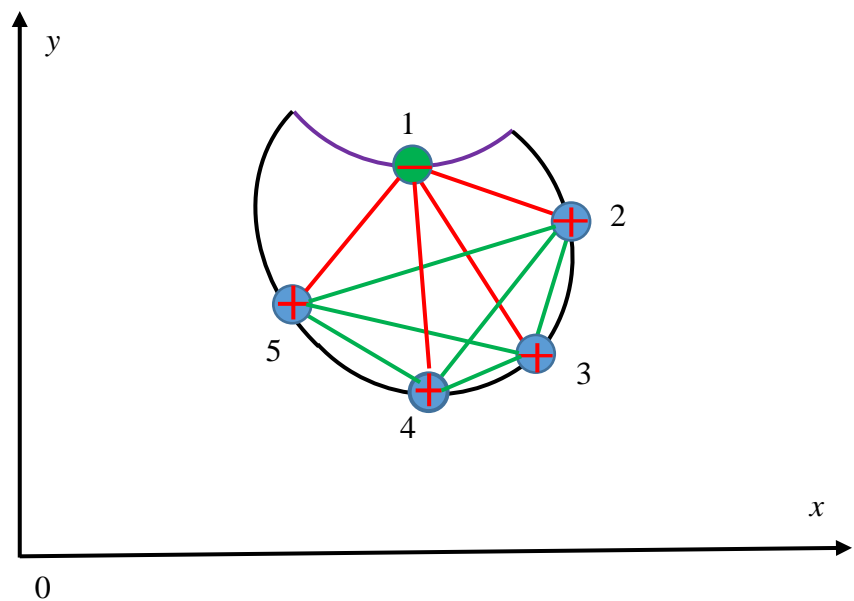

**Figure 3**. A closed curve possessing segments with differently signed curvatures is depicted. The curvature is positive at the black segment of the curvature and negative at the violet segment of the curve. Five points are placed on the curve. Green edges connect the points placed on the segments where the curvatures are of the same sign; red edges connect the points placed on segments were curvatures are differently signed.

Consider first the equally signed vertices. When points numbered "$n$", "$k$" and "$l$" are of the same sign, i.e. if the links $n \leftrightarrow k$ and $k \leftrightarrow l$ are of the same color, link $n \leftrightarrow l$ is necessarily of the same color. Thus, the vertices in the addressed case are connected with the transitive relation (see **Figure 4A**).

Now consider differently signed vertices, When the pairs of vertices numbered "$n$" and "$k$" and "$k$" and "$l$" are differently signed, vertices "$n$" and "$l$" are necessarily of the same sign (see **Figure 4B**). This kind of relations is called in mathematical logic "intransitivity". Let us illustrate this property with the following logical example, involving three groups of experts, labeled "$A$", "$B$" and "$C$" correspondingly. Consider the situation when group of experts "$A$" recognizes group "$B$", and group "$B$" recognizes group "$C$", but group $A$ does not recognize group $C$. In this case, the recognition relation among the expert groups is defined as "intransitive". This is exactly the situation inherent for the relation between three vertices of different signs, shown in **Figure 4B**. It should be emphasized, that no monochromatic triangle will appear when points of various signs of the curvature are located in its vertices; however, the clique built of two monochromic edges will be necessarily present, as shown in **Figure 4B.** On the other hand, the monochromatic triangle will necessarily appear when the three points of the same sign of the curvature are placed in the vertices, as depicted in **Figure 4A**. Using the notions of the Ramsey theory we conclude that $R_{trans,intrans}(2,3) = 3$ is true for the vertices of the graph connected by transitive/intransitive relations.

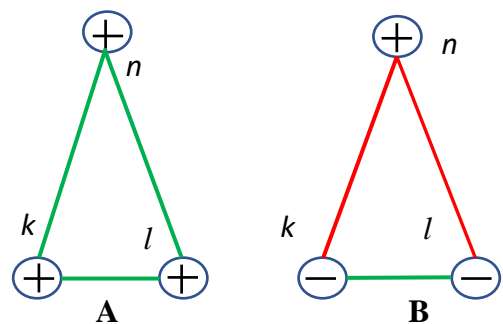

**Figure 4**. Logic interrelations between the points representing positive and negative curvatures are shown. Green edges connect the points placed on the segments where the curvatures are of the same sign; red edges connect the points placed on segments were curvatures are differently signed. **A**. Points in which the curvature is positive are placed in the vertices of the graph. Only green edges appear in the graph. **B**. Points in which the curvature is positive and negative are placed in the vertices of the graph. The monochromatic triangle is impossible in this case.

*2.4 Ramsey approach and the Jordan curves*

Consider one more property of closed curves emerging from the Jordan curve theorem and the Ramsey approach. According to the Jordan theorem every Jordan curve (a plane simple closed curve) divides the plane into an "interior" region bounded by the curve and an "exterior" region containing all of the nearby and far away exterior points. Let us choose six points some of which belong to the interior region and some of which belong to the exterior region as shown in **Figure 5** (the closed curve is shown with the black solid line).

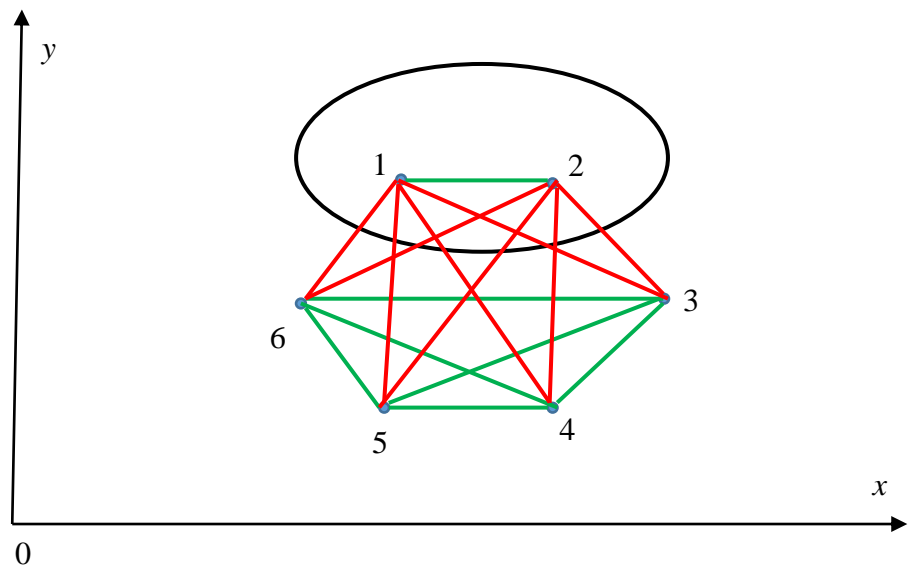

**Figure 5**. Black solid line depicts the Jordan curve dividing the plane into the interior and exterior regions (separated with the solid black line). Points labeled "1" and "2" are located within the interior region, points "3", "4", "5" and "6" are located in the exterior region. Red edges connect points located in different regions; green edges connect points located in the same region. Triangles "456", "345", "356" and "346" are monochromatic.

In the case, illustrated with **Figure 5**, the points numbered "1" and "2" are chosen within the interior region, i.e. within the contour, whereas, the points numbered "2", "4", "5" and "6" are located in the exterior region of the plane, i.e. out of the contour. Let us connect the points located within the same region (whatever interior or exterior) with green edges and the points placed within the different regions with the red eedges, as shown in **Figure 5**. No monochromatic, red triangle is possible in this case. Indeed, if vertices "$n$" and "$k$" and "$k$" and "$l$" are located in different regions, thus, necessarily vertices "n" and "l" are located in the same region. No monochromatic, red triangle is possible in this case, as it is recognized from the **Figure 5**, and we return for the intransitive logical relations for the points located in different regions, which were already discussed in detail in the previous Section. In turn, the logical relation between the points located in the same region are transitive. Thus, monochromatic green triangle will necessarily appear when the three points are placed in the same region of the plane, as it is shown in **Figure 5**. Using the notions of the Ramsey theory we derive $R_{trans,intrans}(2,3) = 3$ for the graph, shown in **Figure 5**.

Now consider three Jordan curves (numbered "*I*", "*II*" and "*III*" in **Figure 6**) dividing the plane in four regions, as depicted in **Figure 6**. Six points labeled "1"… "6" are placed on the plane, as shown in **Figure 6**.

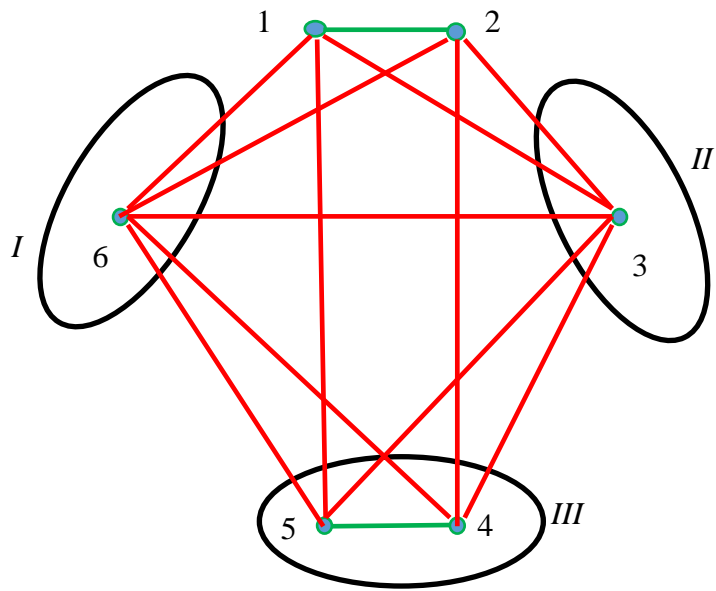

**Figure 6.** Black solid line depicts three closed Jordan curves dividing the plane into the interior and exterior regions Red edges connect points located in different regions; green edges connect points located in the same region. Triangles "136", "236", "346", "356", "234" and "156" are monochromatic.

Let us connect the points located within the same region (whatever its) with green edges and the points placed within the different regions with the red edges, as shown in **Figure 6**. Less than three points are located within the same region (this makes transitive/intransitive relations between vertices impossible for the situation depicted in **Figure 6**). Thus, according to the Ramsey theorem, the monochromatic triangle will necessarily appear in the complete graph, depicted in **Figure 6**. Indeed, triangles "136",

"236", "346", "356", "234" and "156" are monochromatic (red). No monochromatic green triangles are recognized in the graph, shown in **Figure 6**.

Consider one more possible configuration when Jordan curve "*I*" is located within Jordan curve "*II*" as shown in **Figure 7**. In this case, two Jordan curves divide the plane into three distinguishable regions. Consider the set of six points labeled "1"… "6", which are located as shown in **Figure 7**. Again, we connect the points located within the same region (whatever its) with green edges and the points placed within the different regions with the red edges; thus, the complete graph, shown in **Figure 7** emerges.

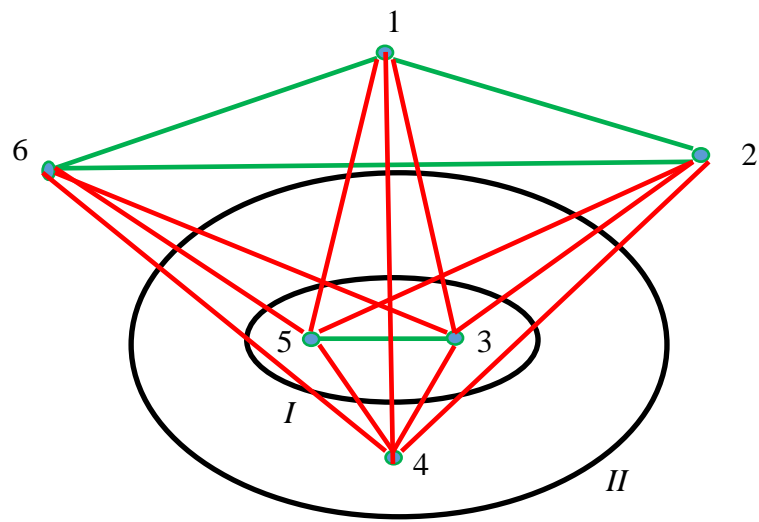

**Figure 7**. The Jordan curve "*I*" is located within Jordan curve "*II*", thus, dividing the plane into three regions. Red edges connect points located in different regions; green edges connect points located in the same region. The triangle "126" is green; triangles "456" and "234" are red.

Three monochromatic triangles, namely green "126" and red "456" and "234" ones, are recognized in the complete graph, shown in **Figure 7**. Green monochromatic triangle "126" should necessarily appear due to the transitivity of logical relations between points "1", "2" and "6", located in the same region; whereas, the red triangles emerge from the Ramsey construction depicted in **Figure 7**.

**Conclusions**

We conclude that the Ramsey theory enables instructive analysis of the discrete sets of points located on closed curves. Monochromatic triangles necessarily appear within the curves, when the points are linked with straight lines. The color of the link is assigned by the slope of the straight line. For set comprising of six points located on the closed curve at least one monochromic triangle should necessarily appear (the Ramsey number

$R(3,3) = 6$). This result immediately provides the non-trivial prediction related to the reflective motion of the point within the dynamical billiard. The second theorem emerges from the combination of the Jordan closed curve and Ramsey theorem. We propose to connect the points located within the same region with green edges and the points placed within the different regions with red edges. In this case, transitivity/intransitivity of the relations between the points should be considered. Ramsey constructions arising from the differential geometry of closed contours are discussed. In this case, also transitive/intransitive Ramsey graphs emerge. Properties of these graphs are discussed.

**CRediT authorship contribution statement**

Nir Shvalb - Writing – review & editing, Writing – original draft, Mark Frenkel - Writing – review & editing, Writing – original draft, Investigation; Formal analysis; Shraga Shoval - Writing – review & editing, Writing – original draft, Supervision; Edward Bormashenko - Writing – review & editing, Writing – original draft; Methodology, Investigation, Formal analysis, Conceptualization.